\documentclass[12pt]{article}
\usepackage{amsmath, amssymb}

\usepackage{setspace}
\begin{document}

\def \AA {\mathcal{A}}
\def \BB {\mathcal{B}}
\def \FF {\mathcal{F}}
\def \HH {\mathcal{H}}
\def \II {\mathcal{I}}
\def \JJ {\mathcal{J}}
\def \KK {\mathcal{K}}
\def \MM {\mathcal{M}}
\def \SS {\mathcal{S}}
\def \VV {\mathrm{v}}
\def \WW {\mathrm{w}}
\def \XX {\mathcal{X}}
\def \CCC {\mathbb{C}}
\def \EEE {\mathbb{E}}
\def \NNN {\mathbb{N}}
\def \PPP {\mathbb{P}}
\def \QQQ {\mathbb{Q}}
\def \RRR {\mathbb{R}}
\def \tr {\mathrm{Tr}}

\def \bk {\mathbf{k}}
\def \bm {\mathbf{m}}
\def \bn {\mathbf{n}}
\def \bp {\mathbf{p}}
\def \bv {\mathbf{v}}

\def \equals {\ = \ }
\def \plus {\ + \ }
\def \minus {\ - \ }

\newtheorem{proposition}{Proposition}

\title{Introduction to determinantal point processes from a quantum probability viewpoint}

\author{Alex D. Gottlieb}

\maketitle

\abstract{Determinantal point processes on a measure space
$(\XX,\Sigma,\mu)$ whose kernels represent trace class Hermitian
operators on $L^2(\XX)$ are associated to ``quasifree" density
operators on the Fock space over $L^2(\XX)$.}

\section{Introduction}
This contribution has been informed and inspired by several
surveys of the topic of determinantal point processes that have
appeared in recent
years.\cite{Soshnikov:Review,Lyons:Review,HoughKrishnapurPeresVirag}
The first of these, Soshnikov (2000), is inspired by the
determinantal point processes that arise in random matrix theory:
the set of eigenvalues of a random matrix is a realization of a
determinantal point process, if the random matrix is sampled from
any of the unitary-invariant ensembles of Hermitian matrices
(e.g., GUE), or from uniform measure on the classical (orthogonal,
unitary, or symplectic) matrix groups, or from the Ginibre
Ensemble. The review by Lyons (2003) is inspired by the Transfer
Current Theorem \cite{BurtonPemantle}, which implies that the
edges occurring in a randomly (uniformly) sampled spanning tree of
a given finite graph $G$ are a determinantal random subset of the
edge set of $G$. Lyons's review concentrates on random subsets of
countable sets, while Soshnikov's review is oriented to treat
discrete subsets of a continuum.  A very recent survey of
determinantal processes (Hough {\it et al.} (2005)) includes the
following newly-found example: the zero set of a power series with
i.i.d. gaussian coefficients is a determinantal point process
\cite{PeresVirag} (the radius of convergence equals $1$ almost
surely).

Hough {\it et al.} (2005) explain how a simple insight gives one a
handle on number fluctuations in determinantal point
processes.\cite{CostinLebowitz,Wie02,Sos00,Sos02}  The insight is
that, in a determinantal point process with finite expected number
of points, the distribution of the number of points is equal to
the distribution of the sum of independent Bernoulli($\lambda_j$)
random variables, where $0 <\lambda_j \le 1$ are the nonzero
eigenvalues of the ``kernel" of the determinantal process.  For
example, consider the number of eigenvalues of a random $n \times
n$ unitary matrix that lie in a given arc $A$ of the unit circle.
Denote this number by $\#_n A$. If the length of $A$ is positive
but less than $2\pi$, then
\begin{equation}
\label{Wieand}
       \frac{\#_n A - \EEE \#_n A }{\sqrt{\ln n}/\pi}
\end{equation}
is asymptotically normal with unit
variance.\cite{Wie02,DiaconisEvans01,Diaconis} The subset of
eigenvalues that lie in $A$ forms a determinantal point process on
$A$, for it is the restriction of a determinantal point process on
the whole circle, hence $\#_n A$ is distributed as a sum of
independent Bernoulli random variables.
  Thus, once one knows that the variance of
$\#_n A$ is $(\ln n)/\pi^2 + o(n)$ \cite{Rains,firstFootnote}, the
asymptotic normality of (\ref{Wieand}) follows from the
Lindeberg-Feller Central Limit Theorem.

Determinantal point processes have a physical interpretation: they
give the joint statistics of noninteracting fermions in a
``quasifree" state. Indeed, this motivated the introduction of the
concept of determinantal (or ``fermion") point processes in the
first place.\cite{Macchi} Analogously defined ``boson" point
processes arise in physics and are called ``permanental" point
processes in probabilistic
writing.\cite{Macchi,HoughKrishnapurPeresVirag} Recently, too,
researchers have continued to investigate determinantal point
fields from a quantum probabilistic point of
view.\cite{Lytvynov,TamuraIto}  We adopt this viewpoint here, and
realize that the satistics of a determinantal point process with
trace class Hermitian kernel $\KK$ on $L^2(\XX)$ are those of
observables on the Fock space $\FF_0(L^2(\XX))$ with respect to
the density operator on $\FF_0(L^2(\XX))$ that determines the
gauge-invariant quasifree state with symbol $\KK$ on the CAR
subalgebra.  However, we do not dwell below on the physical
interpretation, nor do we discuss states on the CAR algebra in the
following.  Our main objective will be to construct the
determinantal point process on $\XX$ with kernel $\KK$, when $\KK$
is the integral kernel of a Hermitian trace class operator on
$L^2(\XX)$ with $0 \le \|K\| \le 1$. Once the construction is
understood, the fact that the number of points in a measurable
subset of $\XX$ is distributed as a sum of independent Bernoulli
random variables becomes obvious.

 Finally, let us remark that determinantal/permanental processes
have a couple of different interesting
generalizations.\cite{DiaconisEvans00,ShiraiTakahashi} And another
rich survey of determinantal processes has just appeared in the
electronic archive!\cite{Johansson}

\section{Determinantal probability measures on finite sets}

Let $\XX$ be a finite set, and let $2^\XX$ denote the set of all
subsets of $\XX$.  Let $\PPP$ denote a probability measure on
$2^\XX$, and let $X$ be a random subset of $\XX$ distributed as
$\PPP$.  Then $\PPP( X \supset E )$ denotes the measure of the
class of all subsets of $\XX$ that contain the subset $E$.  If
there exists a complex-valued function $\KK$ on $\XX \times \XX$
such that
\begin{equation}
     \PPP( X \supset \{x_1,x_2,\ldots ,x_m \} ) \equals \det
     \big(\KK(x_i,x_j)\big)_{i,j=1}^{\ m}
\label{correlations}
\end{equation}
for all subsets $\{x_1,x_2,\ldots ,x_m \}$ of $\XX$, where
$\big(\KK(x_i,x_j)\big)_{i,j=1}^{\ m}$ denotes the $m\times m$
matrix whose $(i,j)^{th}$ entry is $\KK(x_i,x_j)$, then $\PPP$ is
said to be a {\bf determinantal probability measure}
\cite{Lyons:Review} with {\bf kernel} $\KK$.  The probabilities
(\ref{correlations}) determine the probabilities $\PPP(E)$ by
inclusion-exclusion, hence there can be at most one determinantal
probability measure with a given kernel $\KK$.  A very basic
example of a determinantal probability on $2^\XX$ is the law of
the random set produced by independent Bernoulli trials for the
membership of each element of $\XX$; in this case the kernel $
\KK(x',x) \equals  \delta_{x'x} \PPP(x \in X)$.

Suppose $\PPP$ is determinantal with kernel $\KK$.  Then the
complementary probability measure
\[
   \PPP^c(X=S) = \PPP(X=\XX \setminus S)
\]
is determinantal with kernel $\II- \KK$, where $\II(x',x) =
\delta_{x'x}$.   To prove this, use the identity
\begin{eqnarray}
\label{identity}
     &   &
      \det \big(\II(x_i,x_j) - \KK(x_i,x_j)\big)_{i,j=1}^{\ m}  \nonumber \\
      &   &
      \equals 1 \minus \sum_{j=1}^m \KK(x_j,x_j)
      \plus \sum_{1 \le j_1 < j_2 \le m} \det \big( \KK(x_{j_a},x_{j_b}) \big)_{a,b \in \{1,2\}}  \nonumber
      \\
      &   &
      \quad\qquad - \ \sum_{1 \le j_1 < j_2 < j_3 \le m} \det \big( \KK(x_{j_a},x_{j_b}) \big)_{a,b \in \{1,2,3\}}
      \nonumber
      \\
      &   &
       \quad\qquad + \cdots
       \plus (-1)^m \det \big( \KK(x_i,x_j)\big)_{i,j=1}^{\ m} \ .
\end{eqnarray}
The determinants on the right-hand side of (\ref{identity}) are
probabilities according to (\ref{correlations}), therefore $\det
\big(\II(x_i,x_j) - \KK(x_i,x_j)\big)_{i,j=1}^{\ m} $
\begin{eqnarray}
\label{complementaryKernelCalculation}
      & = &
      1 \minus \sum_{j=1}^m \PPP(\{x_j\} \subset X)
      \plus \sum_{1 \le j_1 < j_2 \le m} \PPP(\{x_{j_1},x_{j_2}\} \subset X)   \nonumber
      \\
      &   &
       \plus \cdots \plus (-1)^m \PPP(\{x_1,\ldots,x_m \} \subset X)   \nonumber
      \\
      & = &
      \PPP( X \subset \XX \setminus \{x_1,\ldots,x_m \} )
       \hbox{\qquad\qquad\qquad [by inclusion-exclusion]}\nonumber
      \\
      & = &
      \PPP( (\XX \setminus X) \supset \{x_1,\ldots,x_m \}  )   \nonumber
      \\
      & = &
      \PPP^c( X \supset \{x_1,\ldots,x_m \}  )   \ .
\end{eqnarray}

Suppose that $\XX=\{x_1,x_2,\ldots,x_n\}$ is an $n$-member set.
Define the matrix $K_{ij}=\big( \KK(x_i,x_j)\big)_{ij=1}^{\ n}$.
If $K$ is a Hermitian matrix, then both $K$ and $I-K$ must be
nonnegative matrices, since all of their submatrices have
nonnegative determinants by (\ref{correlations}) and
(\ref{identity},\ref{complementaryKernelCalculation}).  Hence, if
$\KK$ is the kernel of a determinantal random set and $K$ is
Hermitian, then $K$ must be the matrix of a nonnegative
contraction on $\CCC^n$, i.e., necessarily $0 \le \|K\| \le 1$.
Conversely, if $K$ is the matrix of a nonnegative contraction on
$\CCC^n$, then we will show that there exists a determinantal
probability measure on $2^{\{1,\ldots,n\}}$ with kernel
$\KK(i,j)=K_{ij}$.

The rest of this section is devoted to the construction of a
determinantal probability measure whose kernel is a nonnegative
contraction. Our point of view is that there exists a density
operator on the Fock space over $\CCC^n$ whose diagonal elements
in the standard Fock basis give the desired probabilities.

A {\bf density operator} is a nonnegative Hermitian operator of
trace~$1$.

Let $\FF(\CCC^n)$ denote the exterior algebra over $\CCC^n$, i.e.,
\begin{equation}
\label{finiteFock}
   \FF(\CCC^n) \equals \CCC \oplus \CCC^n \oplus \wedge^2 \CCC^n \oplus \cdots \oplus \wedge^{n-1} \CCC^n
      \oplus \wedge^n \CCC^n \ ,
\end{equation}
where $\wedge^m \CCC^n $ denotes the $m^{th}$ exterior power of
$\CCC^n$. The exterior algebra $\FF(\CCC^n)$ is spanned by vectors
of the form $v_1 \wedge v_2 \wedge \ldots \wedge v_m$, where
$v_1,\ldots,v_m$ are any $m$ vectors in $\CCC^n$ and $m$ is any
number between $1$ and $n$ (together with an extra ``vacuum
vector" $\Omega$ to span the first summand). The expression $v_1
\wedge v_2 \wedge \ldots \wedge v_m$ for vectors is formally
multilinear in $v_1,\ldots,v_m$ and satisfies
\[
     v_j \wedge \cdots \wedge v_1 \wedge \cdots \wedge v_m
     \equals - \ v_1 \wedge \cdots \wedge v_j \wedge \cdots \wedge v_m
\]
for $j=2,\ldots,n$. The exterior algebra $\FF(\CCC^n)$ is $2^n$
dimensional and supports the inner product
\[
       \langle v_1 \wedge \cdots \wedge v_{m'} ,\ w_1 \wedge \cdots \wedge w_m \rangle
        \equals
       \delta_{m'm}\ \det \big( \langle v_i, w_j \rangle \big)_{ij=1}^{\ m}
\]
(the vacuum vector is orthogonal to all $v_1 \wedge \ldots \wedge
v_m$ and has unit norm). It can be shown that $\FF(\CCC^n)$ is
isomorphic to a subspace of the {\bf Fock space}
\[
      \FF_0(\CCC^n) \equals \CCC \oplus \CCC^n \oplus (\CCC^n \otimes \CCC^n) \oplus \cdots \oplus (\otimes^n \CCC^n)
\]
via the map that assigns $ 1 \oplus 0_{\CCC^n} \oplus \cdots
\oplus 0_{\otimes^n \CCC^n} $ to $\Omega  $ and
\begin{equation}
\label{isomorphism}
     0_{\CCC} \oplus \cdots \oplus 0_{\otimes^{m-1} \CCC^n} \oplus
            \mathcal{S}\ell [v_1,\ldots,v_m] \oplus 0_{\otimes^{m+1} \CCC^n} \oplus \cdots \oplus 0_{\otimes^n \CCC^n}
\end{equation}
to $v_1 \wedge v_2 \wedge \ldots \wedge v_m$.  In
(\ref{isomorphism}), $\mathcal{S}\ell [v_1,\ldots,v_m]$ denotes
the Slater determinant
\begin{equation}
\label{Slater}
      \mathcal{S}\ell [v_1,\ldots,v_n] \equals
      \frac{1}{\sqrt{n!} }
      \sum_{\pi \in \SS_m} \hbox{sgn}(\pi)\ U_{\pi}(v_1 \otimes v_2 \otimes \cdots \otimes v_m)\ ,
\end{equation}
where $\SS_m$ denotes the group of permutations of
$\{1,\ldots,m\}$ and $U_{\pi}$ is the unitary operator defined on
$\otimes^m \CCC^n$ when $\pi \in \SS_m$ by the condition that
\begin{equation}
\label{perm}
    U_{\pi}( w_1 \otimes w_2 \otimes \cdots \otimes w_m ) \equals \ w_{\pi^{-1}(1)} \otimes w_{\pi^{-1}(2)} \otimes \cdots \otimes w_{\pi^{-1}(m)}
\end{equation}
for all $w_1,\ldots,w_m \in \CCC^n$. Henceforth, we identify the
exterior algebra $\FF(\CCC^n)$ with this subspace of
$\FF_0(\CCC^n)$, and call it the ``fermion Fock space."

An orthonormal basis of $\FF(\CCC^n)$, called a {\bf Fock basis}
or ``occupation number" basis, can be built using any ordered
orthonormal basis $\VV = (v_1,\ldots,v_n)$ of $\CCC^n$. The
vectors of the Fock basis can be conveniently indexed by subsets
of $\{1,\ldots, n\}$: the empty subset of $\{1,\ldots, n\}$
corresponds to the vacuum vector $\Omega$ and a nonempty subset
$\{j_1,\ldots,j_m \} \subset \{1,\ldots, n\}$ with $j_1 < \ldots <
j_m$ corresponds to the vector $ v_{j_1} \wedge  \cdots \wedge
v_{j_m} $.  That is, the orthonormal set $\{f_{\VV}(S) \ | \ S
\subset \{1,\ldots,n\}\}$ is a basis for $\FF(\CCC^n)$, where $
f_{\VV}(\{\}) = \Omega $ and $     f_{\VV}(S) = v_{j_1} \wedge
\cdots \wedge v_{j_m} $ when $S=\{j_1,\ldots,j_m\}$ with $ j_1 <
\ldots < j_m$.

Suppose $K$ is a nonnegative contraction on $\CCC^n$ and let $\VV
= (v_1, \ldots, v_n )$ be an ordered orthonormal basis of $\CCC^n$
such that $K v_j = \lambda_j v_j$ for all $j$. Let $D_K$ denote
the density operator
\begin{equation}
\label{FockDensity}
     D_K \equals \sum_{S \subset \{1,\ldots,n\}}
     \Big\{  \prod_{ k \in S }\lambda_k \prod_{ k \notin S }(1-\lambda_k) \Big\}
     \ \big\langle f_{\VV}(S), \ \cdot\ \big\rangle
     f_{\VV}(S)
\end{equation}
on $\FF(\CCC^n)$, where $\big\langle f_{\VV}(S), \ \cdot\
\big\rangle f_{\VV}(S)$ denotes the rank-one orthogonal projector
onto the span of $f_{\VV}(S)$.

\begin{proposition}
Let $K$ be a nonnegative contraction on $\CCC^n$ and let $D_K$
denote the associated density operator (\ref{FockDensity}) on the
Fock space $\FF(\CCC^n)$. Then, for all ordered orthonormal bases
$\WW=(w_1,\ldots,w_n)$ of $\CCC^n$,
\begin{equation}
\label{probability}
      S \ \longmapsto \ \big\langle f_{\WW}(S), D_K f_{\WW}(S) \big\rangle
\end{equation}
is a determinantal probability measure on $2^{\{1,\ldots,n\}}$
with kernel $\KK(i,j) = \langle K w_i,  w_j \rangle $.
\end{proposition}
\noindent {\bf Proof:} \qquad We first define the ``second
quantization" maps from operators $A$ on $\otimes^m\CCC^n$ to
operators $\Gamma_m[A]$ on the Fock space $\FF_0(\CCC^n)$, and the
dual maps from density operators $D$ on $\FF_0(\CCC^n)$ to
$m$-particle ``correlation operators" $\mathrm{K}_m[D]$ on
$\otimes^n\CCC^n$.

Let $\JJ(m,k)$ denote the set of injections of $\{1,\ldots,m\}$
into $\{1,\ldots,k\}$.  The cardinality of $\JJ(m,k)$ is $
k^{[m]} \equiv k(k-1)\cdots (k-m+1)  $, the $m^{th}$ factorial
power of $k$. For any operator $A$ on $\otimes^m \CCC^n$, and any
injection $j \in \JJ(m,k)$ with $k \ge m$, we define the operator
\[
     A^{(j)} \equals
     U_{(1 j_1)(2 j_2)\cdots(m j_m)} (A \otimes I \otimes \cdots \otimes I ) U_{(1 j_1)(2 j_2)\cdots(m j_m)}
\]
on $\otimes^k \CCC^n$, where $U_{(1 j_1)(2 j_2)\cdots(m j_m)}$
denotes the permutation operator (\ref{perm}) for the product of
disjoint transpositions $(1 j_1)(2 j_2)\cdots (m j_m)$. Define
$\Gamma_m[A]$ on $\FF_0(\CCC^n)$ by
\[
     \Gamma_m[A] \equals   0_{\CCC} \oplus \cdots \oplus 0_{\otimes^{m-1} \CCC^n}
     \oplus \sum_{j \in \JJ(m,m)} A^{(j)} \oplus  \cdots \oplus  \sum_{j \in \JJ(m,n)} A^{(j)}   \  .
\]
A density operator $D$ on $\FF(\CCC^n)$ extends to a density
operator $D \oplus 0$ on $\FF_0(\CCC^n)$, which we will denote by
$D$ as well.  The map $A \longmapsto \tr( D \Gamma_m[A])$ is a
linear functional on the space of linear operators on $\otimes^m
\CCC^n$. Therefore there exists a unique operator
$\mathrm{K}_m[D]$ on $\otimes^m \CCC^n$ such that
\[
\tr( \Gamma_m[A] D) = \tr( A \mathrm{K}_m[D])
\]
for all linear operators $A$ on $\otimes^m \CCC^n$. In physics
language, $\mathrm{K}_m[D]$ is the $m$-particle correlation
operator for the state with density operator $D$.  If $D_K$ is
defined as in (\ref{FockDensity}), the key identity
\begin{equation}
\label{key}
    \mathrm{K}_m[D_K] \equals  \stackrel{ m \mathrm{\ times} }{\overbrace{K \otimes K \otimes \cdots \otimes K} }
    \sum_{\pi \in \SS_m} \hbox{sgn}(\pi)U_{\pi}
\end{equation}
may verified by comparing matrix elements of both sides with
respect to the basis $\{ v_{j_1}\otimes \cdots \otimes v_{j_m} |\
j_1,\ldots,j_m \in \{1,\ldots,n\}\}$.

Given an ordered orthonormal basis $\WW = (w_1,\ldots,w_n)$, let
$P^{\WW}_j$ denote the projector $\langle w_j ,\cdot \rangle w_j $
for $j=1,\ldots,n$. For distinct $ x_1,\ldots, x_m $, the operator
$\Gamma_m[P^{\WW}_{x_1}\otimes \cdots \otimes P^{\WW}_{x_m}]$ is
diagonal in the Fock basis $\{f_{\WW}(S) \ | \ S \subset
\{1,\ldots,n\}\}$, and
\[
    \Gamma_m[P^{\WW}_{x_1}\otimes \cdots \otimes P^{\WW}_{x_m}]f_{\WW}(S) \equals
     \Bigg\{ \begin{array}{cl}
          f_{\WW}(S)   & \hbox{if } \{x_1,\ldots,x_m\} \subset S \\
          0   & \hbox{otherwise. }  \\
         \end{array}
\]
Let $\PPP(X=S)$ denote the probability (\ref{probability}).   Then
\begin{eqnarray*}
\PPP( X \supset \{x_1,x_2,\ldots ,x_m \} ) & = &   \sum_{S \supset
\{x_1,x_2,\ldots ,x_m \}} \PPP( X = S  )
\\
& = &
\sum_{S \supset \{x_1,x_2,\ldots ,x_m \}} \big\langle f_{\WW}(S), D_K f_{\WW}(S) \big\rangle  \\
& = & \tr( \Gamma_m[P^{\WW}_{x_1}\otimes \cdots \otimes
P^{\WW}_{x_m}] D_K)
\\
& = &
\tr(  (P^{\WW}_{x_1}\otimes \cdots \otimes P^{\WW}_{x_m}) \mathrm{K}_m[D_K] ) \\
& = &
\tr \Big(  (P^{\WW}_{x_1}K \otimes \cdots \otimes P^{\WW}_{x_m}K  )  \sum_{\pi \in \SS_m} \hbox{sgn}(\pi)U_{\pi}  \Big) \\
& = & \det \big( \langle K w_{x_i},  w_{x_j} \rangle
\big)_{i,j=1}^{\ m} \equals \det \big( \KK(x_i,x_j)
\big)_{i,j=1}^{\ m} \ .
\end{eqnarray*}
This proves the proposition.  \hfill $\square$

\section{Determinantal finite point processes}

A {\bf finite point process on} $\XX$ is a random finite subset of
a space $\XX$. Let $\Sigma$ be a $\sigma$-field of measurable
subsets of $\XX$.\cite{footnote}
 A finite point process on $(\XX,\Sigma)$ is specified by the probabilities
 $p_0,p_1,p_2,\ldots$ that there are $0,1,2,\ldots$ points in the configuration,
 and, for each $n$ such that $p_n\ne 0$, a symmetrical conditional probability measure
 $\rho_n$ on $ (\XX^n, \otimes^n \Sigma)$.\cite{DaleyVere-Jones}
Now let $\mu$ be any positive ``reference" measure on
$(\XX,\Sigma)$. A finite point process is {\bf determinantal on}
$(\XX,\Sigma,\mu)$ with kernel $\KK:\XX\times \XX \longrightarrow
\CCC$ if
\begin{equation}
\label{correlationsPointProcess}
    \EEE\Big( \prod_{j=1}^m \#(X \cap E_j) \Big) \equals
    \int_{E_1}\cdots\int_{E_m}  \det \big( \KK(x_i,x_j ) \big)_{ij=1}^{\ m} \mu(dx_1) \cdots \mu(dx_m)
\end{equation}
for all disjoint, measurable $E_1,\ldots,E_m$, $m\ge
1$.\cite{HoughKrishnapurPeresVirag} If $\KK(x,y)$ is the standard
version of the integral kernel of a nonnegative trace class
contraction $K$ on $L^2(\XX,\Sigma,\mu)$, then there exists a
unique \cite{footnoteUniqueness} determinantal point process on
$(\XX,\Sigma,\mu)$ with kernel $\KK(x,y)$.   Conversely, if the
kernel of a determinantal point process on $\XX$ is the integral
kernel of a trace class Hermitian operator $K$ on $L^2(\XX)$, then
$K$ must be a nonnegative contraction.\cite{footnoteConverse}

In this section we construct the determinantal point process on
$\XX$ whose kernel is the standard kernel (\ref{canonical}) of a
given trace class operator $K$ on $L^2(\XX)$ with $0 \le \|K\| \le
1$. This is accomplished by constructing a density operator on
$\FF(L^2(\XX))$ as we have done in the preceding section --- our
quantum probabilistic point of view. There are many other ways to
accomplish the same end, with or without our point of view. The
original approach of Macchi (1975) was to start with a formula for
the Janossy densities \cite{JanossyDensities} of the desired point
process, and then to verify (\ref{correlationsPointProcess}) for
that process. Soshnikov (2000) attacks the problem by first
showing that certain Fredholm determinants involving $K$ define
factorial moment generating functions for finite families of
random variables $\{\#(X \cap E_j)| E_j \in \Sigma\}$, and then
constructing the desired determinantal point process via
Kolmogorov extension from its finite dimensional distributions.
Lyons (2003) uses the geometry of Fock space, but only in the case
where $K$ is a finite rank projector, then dilates nonnegative
contractions to projections on a larger space to handle the
general case.  Hough {\it et al.} (2005) verify directly that
kernels of finite rank projectors yield determinantal point
processes, then treat the general case as a mixture, in the
probabilistic sense, of determinantal processes with projector
kernels.

Any density operator on $\FF(L^2(\XX,\Sigma,\mu))$ of the form
$D=\oplus D_n$ defines a finite point process on $\XX$ as follows.
 $p_n=\tr(D_n)$ is the probability of the event the configuration has exactly $n$ points.
The measure $\rho_n$ is absolutely continuous with respect to
$\otimes^n \mu$ and it is defined by way of the isomorphism
\[
    \otimes^n L^2(\XX,\Sigma,\mu) \ \stackrel{\sim}{=}\ L^2(\XX^n,\otimes^n\Sigma, \otimes^n\mu).
\]
Regarding $D_n$ as an operator on $L^2(\XX^n,\otimes^n\Sigma,
\otimes^n\mu)$, define
 $\rho_n(E) = p_n^{-1} \tr(D_n \MM_E )$,
where $ \MM_E $ denotes the operator on $L^2(\XX^n)$ of
multiplication by the indicator function of $E \in \otimes^n
\Sigma$.

Given a trace class nonnegative contraction $K$ on $L^2(\XX)$, a
density operator $D_K$ on $\FF(L^2(\XX))$ may be defined using the
spectral information in $K$ as was done in (\ref{FockDensity})
above:
\begin{equation}
     D_K \equals \sum_{n=0}^\infty \sum_{\stackrel{S \subset \{1,2,\ldots\}:}{ \# S =n}}
     \Big\{ \prod_{ k \in S }\lambda_k \prod_{ k \notin S }(1-\lambda_k) \Big\}
     \ \big\langle f(S), \ \cdot\ \big\rangle     f(S)\ ,
\label{thereItIs}
\end{equation}
where $\lambda_1 \ge \lambda_2 \ge \ldots $ are the eigenvalues of
$K$ and $\{f(S)\}$ is the Fock basis constructed from the
eigenvectors of $K$ (really, any extension of an orthonormal
system of eigenvectors of $K$ to an orthonormal basis of
$L^2(\XX)$). $D_K$ has the form $\oplus (D_{K})_n$ and it can be
verified \cite{self} that the $m$-particle correlation operator
$\mathrm{K}_m[D_K]$ exists and satisfies the key identity
(\ref{key}). Let $p_n$ and $\rho_n$ be as defined above, for
$D_K$. Let $\EEE$ denote the expectation with respect to random
point process defined by these $p_n$ and $\rho_n$.  Then
\begin{equation}
\label{first}
    \EEE\Big(  \prod_{j=1}^m \#(X \cap E_j) \Big)
    \equals
       \tr( \Gamma_m[ P_1 \otimes \cdots \otimes P_m ] D_K )  \ ,
\end{equation}
where $P_j$ denotes the orthogonal projector $\MM_{E_j}$ on
$L^2(\XX)$, for both sides of (\ref{first}) equal
  $ \sum\limits_{n \ge m} p_n \int_{\XX^n}   \prod_{j=1}^m \#(\{x_1,\ldots,x_n\} \cap E_j)   \rho_n(dx_1 \cdots dx_n) $.
But
\begin{equation}
\label{second}
    \tr( \Gamma_m[ P_1 \otimes \cdots \otimes P_m ] D_K )
    \equals
      \tr\Big(  (P_1 K \otimes \cdots \otimes P_m K) \sum_{\pi \in \SS_m} \hbox{sgn}(\pi)U_{\pi}  \Big)
\ ,
\end{equation}
since $\mathrm{K}_m[D_K ]=\big(\otimes^mK \big) \sum\limits_{\pi
\in \SS_m} \hbox{sgn}(\pi)U_{\pi} $. Finally, one may verify
\cite{anotherFootnote} that
\begin{equation}
\label{third}
    \tr\Big(  \big( \otimes_{j=1}^m P_j  K \big) \sum_{\pi \in \SS_m} \hbox{sgn}(\pi)U_{\pi}  \Big)
    \equals
       \int_{E_1}\cdots\int_{E_m}  \det \big( \KK(x_i,x_j ) \big)_{ij=1}^{\ m} \mu(dx_1) \cdots \mu(dx_m)
\end{equation}
if $\KK(x,y)$ is the usual version of the integral kernel of $K$,
i.e.,
\begin{equation}
\label{canonical}
    \KK(x,y) \equals \sum_j \lambda_j \phi_j(x) \overline{\phi_j(y)}
\end{equation}
where $K\phi_j = \lambda_j \phi_j$ and $\sum \lambda_j = \tr K$.
Equations (\ref{first}) - (\ref{third}) imply
(\ref{correlationsPointProcess}) holds; the finite point process
defined through $D_K$ is determinantal with kernel $\KK$.

Now that we have constructed the process, we can see immediately
from (\ref{thereItIs}) that the total number of points in a random
configuration is distributed as the sum of Bernoulli($\lambda_j$)
random variables.  In particular, $\prod (1-\lambda_k)$ is the
probability that there are no points at all.  This equals the
Fredholm determinant $\hbox{Det} (I-K)$. Let $E$ be a measurable
subset of $\XX$.  It is not difficult to check via
(\ref{correlationsPointProcess}) that the determinantal point
process on $(E,\Sigma|_E,\mu|_E)$ with kernel $K_E \equiv \MM_C K
\MM_C \big|_E$ is the restriction to $E$ of the determinantal
point process with kernel $K$ on $\XX$. Hence the probability that
there are no points in $E$ equals the Fredholm determinant
$\hbox{Det}(I-K_E)$. In the context of random matrix theory, this
yields formulas for the spacing distributions of
eigenvalues.\cite{TracyWidom93,TracyWidom98}

In case $\|K\|<1$, set $L=(I-K)^{-1} K$.  It is easy to check from
(\ref{thereItIs}) that
\[
     (D_K)_n \equals  \frac{1}{n!}\ \hbox{Det}(I-K) \big\{ \otimes^n L \big\}
      \sum_{\pi \in \SS_n} \hbox{sgn}(\pi)U_{\pi}
\]
by comparing matrix elements of both sides of this identity with
respect to an eigenbasis of $K$. This identity yields the
determinantal formulas for the Janossy
densities.\cite{Macchi,DaleyVere-Jones}

\section{Determinantal processes of infinitely many points}
Suppose that $\XX$ is a locally compact Hausdorff space satisfying
the second axiom of countability
\cite{SecondCountableLocallyCompactHausdorff}, and let $\Sigma$
denote the Borel field of $\XX$. In this context, a point process
is a random nonnegative integer-valued Radon measure (a Radon
measure is a Borel measure which is finite on any compact
set).\cite{HoughKrishnapurPeresVirag}  Let $\mu$ be a
$\sigma$-finite Radon measure on $\XX$.\cite{Sos02} A point
process on $(\XX,\Sigma,\mu)$ is {\bf determinantal} with kernel
$\KK$ if (\ref{correlationsPointProcess}) holds.

Most work on determinantal processes with infinitely many points
has been done for the cases where $\XX$ is a countable set with
the discrete topology and $\mu$ is counting measure, or $\XX$ is a
connected open subset of $\RRR^d$ and $\mu$ is Lebesgue measure,
or $\XX$ is a finite disjoint union or Cartesian product of said
spaces. If $K$ is a {\it locally} trace class Hermitian operator
on $L^2(\XX)$ such that $0 \le \|K\| \le 1$, then (a version of)
its integral kernel is the kernel of a determinantal point process
on $\XX$.\cite{Soshnikov:Review}  This point process is the limit
in distribution of the determinantal processes with kernels ${\bf
1}_C(x) \KK(x,y)  {\bf 1}_C(y)$, where $C$ ranges over an
increasing family of compact subsets of $\XX$. Conversely, if the
kernel of a locally trace class Hermitian operator $K$ defines a
determinantal point process, then $K$ must be a nonnegative
contraction.\cite{Soshnikov:Review,Macchi}

The case of a countably infinite set $\XX$ with counting measure
is treated in detail in Lyons (2003).  In this pleasant special
case, Hermitian operators on $L^2(\XX,2^{\XX},\#)$ are
automatically locally trace class.  On the other hand, equations
(\ref{correlations}) and
(\ref{identity})-(\ref{complementaryKernelCalculation}) readily
imply that $K$ and $I-K$ are both nonnegative operators.
Therefore, the kernel of a Hermitian operator $K$ on
$L^2(\XX,2^{\XX},\#)$ is the kernel of a determinantal point
process on $(\XX,2^{\XX},\#)$ if and only if $K$ is a nonnegative
contraction.

\section*{Acknowledgments}
This work was supported by the Austrian Ministry of Science
(BM:BWK) via its grant for the Wolfgang Pauli Institute and by the
Austrian Science Foundation (FWF) via the START Project
(Y-137-TEC) of N.J. Mauser. I dedicate this to Steve Evans who
introduced me to determinantal point processes.


\begin{thebibliography}{0}
\bibitem{footnoteConverse}
The hypothesis that $K$ is trace class is used to prove that
$\|K\| \le 1$ in
general.\cite{Soshnikov:Review,HoughKrishnapurPeresVirag} However,
this hypothesis is not required if $\XX$ is a locally compact
space and $\KK$ is continuous (and $\mu$ is a Radon measure on the
Borel field).  In such cases, if the (continuous) kernel of a
determinantal point process on $\XX$ is the integral kernel of a
bounded --- but not necessarily trace class --- Hermitian operator
$K$ on $L^2(\XX)$, then necessarily $0 \le \|K\| \le 1$.
\bibitem{JanossyDensities}  The Janossy densities are the densities of the measures $n!p_n\rho_n$.
\bibitem{SecondCountableLocallyCompactHausdorff}  Such are the spaces considered in the work of A. Lenard on
correlation densities for infinitely many particles (viz.
Theorem~1 in \cite{Soshnikov:Review}).
\end{thebibliography}
\end{document}